\documentclass[12pt]{article}
\usepackage[margin=2cm]{geometry}
\usepackage{amsfonts}
\usepackage{amssymb}
\usepackage{amsthm}
\usepackage{amsmath}
\usepackage{latexsym}
\usepackage{color}
\usepackage{mathrsfs}

\begin{document}
\allowdisplaybreaks[4]
\newtheorem{thme}{Theorem}
\newtheorem{lemma}{Lemma}
\newtheorem{pron}{Proposition}
\newtheorem{re}{Remark}
\newtheorem{thm}{Theorem}
\newtheorem{Corol}{Corollary}
\newtheorem{exam}{Example}
\newtheorem{defin}{Definition}
\newtheorem{remark}{Remark}
\newtheorem{property}{Property}
\newcommand{\bco}{\color{blue}}
\newcommand{\rco}{\color{red}}
\newcommand{\la}{\frac{1}{\lambda}}
\newcommand{\sectemul}{\arabic{section}}
\renewcommand{\thethme}{\sectemul.\Alph{thme}}
\renewcommand{\theequation}{\sectemul.\arabic{equation}}
\renewcommand{\thepron}{\sectemul.\arabic{pron}}
\renewcommand{\thelemma}{\sectemul.\arabic{lemma}}
\renewcommand{\there}{\sectemul.\arabic{re}}
\renewcommand{\thethm}{\sectemul.\arabic{thm}}
\renewcommand{\theCorol}{\sectemul.\arabic{Corol}}
\renewcommand{\theexam}{\sectemul.\arabic{exam}}
\renewcommand{\thedefin}{\sectemul.\arabic{defin}}
\renewcommand{\theremark}{\sectemul.\arabic{remark}}
\def\REF#1{\par\hangindent\parindent\indent\llap{#1\enspace}\ignorespaces}
\def\lo{\left}
\def\ro{\right}
\def\be{\begin{equation}}
\def\ee{\end{equation}}
\def\beq{\begin{eqnarray*}}
\def\eeq{\end{eqnarray*}}
\def\bea{\begin{eqnarray}}
\def\eea{\end{eqnarray}}
\def\r{random walk}
\def\o{\overline}

\title{\large\bf A necessary and sufficient condition for the subexponentiality of product convolution
}
\author{\small Hui Xu~~ Fengyang Cheng~~ Yuebao Wang\thanks{Research supported by the National Natural Science Foundation of China (Nos. 11401415, 11071182) and Jiangsu Overseas Research and Training Program for Prominent University Young and Middle-aged Teachers and Presidents}
\thanks{Corresponding author.
Telephone: 86 512 67422726. Fax: 86 512 65112637. E-mail:
ybwang@suda.edu.cn}~~Dongya Cheng\\
{\footnotesize School of Mathematical Sciences, Soochow
University, Suzhou 215006, China}}
\date{}
\maketitle

\begin{center}
{\noindent\small {\bf Abstract }}
\end{center}

{\small
Let $X$ and $Y$ be two independent and nonnegative random variables with corresponding distributions $F$ and $G$.
Denote by $H$ the distribution of the product  $XY$, called the product convolution of $F$ and $G$.
Cline and Samorodnitsky (1994) proposed sufficient conditions for $H$ to be subexponential, given the subexponentiality of $F$.
Relying on a related result of Tang (2008) on the long-tail of product convolution,
we obtain a necessary and sufficient condition for the subexponentiality of $H$, given that of $F$.
We also study the reverse problem and obtain sufficient conditions for the subexponentiality of $F$ given that of $H$.
Finally, we apply the obtained results to the asymptotic study of the ruin probability in a discrete-time insurance risk model with stochastic returns.
\medskip

{\it Keywords:} necessary and sufficient condition; product convolution;
subexponential distribution; discrete-time insurance risk model

\medskip

{\it AMS 2010 Subject Classification:} Primary 60E05, secondary 62E20, 60G50.}

\section{Introduction and main results}
\setcounter{thm}{0}\setcounter{Corol}{0}\setcounter{lemma}{0}\setcounter{pron}{0}\setcounter{equation}{0}
\setcounter{re}{0}\setcounter{exam}{0}\setcounter{property}{0}\setcounter{defin}{0}

Throughout the paper, we assume that $X$ and $Y$ are two independent and nonnegative random variables with corresponding
distributions $F$ and $G$.
Denote by $H$ the distribution of the product $XY$, called the product convolution of $F$ and $G$,
as opposed to the usual sum convolution $F*G$ defined to be the distribution of the sum $X+Y$.

The concept of product convolution plays an important role in applied probability and has become increasingly interesting in recent years; see Cline and Samorodnitsky (1994), Tang (2006), Li and Tang (2015), Samorodnitsky and Sun (2016), among many others. Note that the present value of a random future claim, which is one of most fundamental quantities in finance and insurance, is expressed as the product of the random claim amount and the corresponding stochastic present value factor. Thus, the study of the tail behavior of product convolutions has immediate implications in finance and insurance. In this paper, we are interested in the closure property of product convolution within the subexponential class. More precisely, given that $F$ is subexponential, we pursue a necessary and sufficient condition such that $H$ is subexponential too.

Throughout the paper, the following notation and conventions are in force. For a distribution $V$,
denote by $\overline{V} = 1-V$ its tail distribution.
Unless otherwise stated, all limiting relations are according to $x\to\infty$. For two positive functions
$f(\cdot)$ and $g(\cdot)$, the notation $f(x)\sim g(x)$ means that $\lim f(x)/g(x)=1$.
With $s=\limsup f(x)/g(x)$, we write $f(x)=o\big(g(x)\big)$ if $s=0$, $f(x)=O\big(g(x)\big)$ if $s<\infty$, and $f(x)\lesssim g(x)$ if $s\leq 1$.
Finally, we write $f(x)\asymp g(x)$ if $f(x)=O\big(g(x)\big)$ and $g(x)=O\big(f(x)\big)$.

A distribution $V$ supported on $(-\infty,\infty)$ is said to belong to the class $\mathcal{L}(\gamma)$ for some $\gamma\ge0$ if, for any $t$,
\begin{eqnarray*}
\overline V(x-t)\sim e^{\gamma t}\overline V(x).
\end{eqnarray*}
When $\gamma>0$ and the distribution $V$ is lattice, the variables $x$ and $t$ above should be restricted to values of the lattice span. When $\gamma=0$, it reduces to $\mathcal{L}=\mathcal{L}(0)$, the class of long-tailed distributions. It is well known that every $V\in\mathcal{L}$ is heavy tailed in the sense
of infinite exponential moments.

A distribution $V$ supported on $[0,\infty)$ is said to belong to the subexponential class $\mathcal{S}$ if
\begin{eqnarray*}
\overline{V^{*2}}(x)=\overline{V*V}(x)\sim2\overline V(x).
\end{eqnarray*}
Furthermore, a distribution $V$ supported on $(-\infty,\infty)$ is still said to be subexponential if the distribution $V_+$ defined by
$$V_+(x)=V(x)\textbf{1}(x\ge0)$$
is subexponential, where $\textbf{1}(A)$ denotes the indicator  of an event $A$, equal to $1$ if $A$ occurs and  to $0$ otherwise.
It is well known that any subexponential distribution is long tailed.
The class $\mathcal{S}$ was introduced by Chistyakov (1964); see Embrechts et al. (1997), Asmussen (2000) and Foss et al. (2013) for extensive discussions on this class and its applications.

Another well-known class of heavy-tailed distributions is the class $\mathcal{D}$ of dominated variation introduced by Feller (1971). Recall that a distribution $V$ supported on $(-\infty,\infty)$ belongs to the class $\mathcal{D}$ if
$$\overline {V}(tx)=O\big(\overline V(x)\big)$$ %.
holds for some (and hence for all) $0<t<1$.

The intersection of the classes $\mathcal{L}$ and $\mathcal{D}$ contains the famous class $\mathcal{R}$ of distributions of regular variation.
By definition, a distribution $V$ is said to belong to the class $\mathcal{R}$ if its tail $\overline V$ is regularly varying, namely, if
for some $\alpha\ge0$ and all $t>1$,
\begin{eqnarray*}
\overline {V}(tx)\sim {t}^{-\alpha}\overline V(x).
\end{eqnarray*}
The reader is referred to Bingham et al. (1987) and Resnick (1987) for textbook treatments of functions of regular variation.

Moreover, Tang (2006) introduced the class $\cal{A}$ consisting of all subexponential distributions $V$ that satisfy the following property: for some $t>1$,
$$\limsup \overline V(tx)/\overline V(x)<1.$$
As Tang (2006) pointed out, the class $\cal{A}$ is just slightly smaller than the class $\cal{S}$ and it essentially does not exclude any useful subexponential distributions.

Now recall several works on closure properties of the product convolution. Embrechts and Goldie (1980) showed that if $F\in\mathcal{R}$ and
$\overline G(x)=o\big(\overline F(x)\big)$, then $H\in\cal{R}$. Cline and Samorodnitsky (1994)
studied the closure property of the product convolution in the class $S$. Given $F\in\mathcal{S}$,
they established certain sufficient conditions for $H\in\cal{S}$.
Samorodnitsky and Sun (2016) extended the study to multivariate subexponential distributions. Tang (2006) investigated the closure property of the product convolution in the class $\cal{A}$ under some more relaxed sufficient conditions. Tang (2008) presented a necessary and sufficient condition for $H\in\cal{L}$ when $F\in\cal{L}(\gamma)$ for some $\gamma\ge0$.

Now we recall Theorem 2.1 of Cline and Samorodnitsky (1994) in detail:
\vspace{0.3cm}
\begin{thme}\label{1a} (Cline and Samorodnitsky (1994)).\\
Assume that $F\in\mathcal{S}$ and there is a function $a(\cdot):[0,\infty)\mapsto(0,\infty)$ satisfying:

(a) $a(x)\uparrow\infty$;

(b) $a(x)/x\downarrow 0$;

(c) $\overline F\big(x-a(x)\big)\sim \overline F(x)$;

(d) $\overline G\big(a(x)\big)=o\big(\overline H(x)\big)$.

\noindent Then $H\in\mathcal{S}$.
\end{thme}

Here are some notes on Theorem \ref{1a}. First, Cline and Samorodnitsky (1994) pointed out that, for a distribution $F\in\mathcal{L}$,
there always exists a function $a(\cdot):[0,\infty)\to(0,\infty)$ satisfying conditions (a)-(c), namely, the tail $\overline F$ is `$a$-insensitive'.
Second, if $G$ has a bounded support, then condition (d) becomes trivial while conditions (a), (b) and (c) are automatically implied by $F\in\mathcal{S}$. For this case, Theorem \ref{1a} reduces to Corollary 2.5 of Cline and Samorodnitsky (1994). Third, Tang (2006) noted that conditions (a), (b) and (d) are satisfied if and only if
\begin{equation}\label{blxtj2}
\overline G(x)=o\big(\overline H(bx)\big)\ \text{for every}\ b>0.
\end{equation}
We will construct an example in Subsection 3.1 in which $F\in\mathcal{R}$, $G\notin\mathcal{L}$ and $H\in\mathcal{S}$, while relation (\ref{blxtj2}) does not hold. This example indicates that condition (d) in Theorem 1.A is not an necessary condition for $H\in\mathcal{S}$.

A question naturally arises: When $F\in\mathcal{S}$ and $G$ has an unbounded support, can we establish a sufficient and necessary condition for $H\in\mathcal{S}$? Now we recall Theorem 1.1 of Tang (2008). Denote by $D[V]$ the set of all positive points of discontinuity of the distribution $V$.
\begin{thme}\label{1b} (Tang (2008))\\
Assume that $F\in\mathcal{L}(\gamma)$ for some $\gamma\ge0$, and $G$ has an unbounded support $[0,\infty)$.
Then $H\in\mathcal{L}$ if and only if one of the following holds: either $D[F]=\emptyset$, or $D[F]\not=\emptyset$ and
\begin{equation}\label{blxtj}
\overline G(x/d)-\overline G\big((x+1)/d\big)=o\big(\overline H(x)\big)\ \text{for all}\ d\in D[F].
\end{equation}
\end{thme}

It turns out that the condition used in Theorem \ref{1b} is the sufficient and necessary condition in the following our main result with a different premise.

\begin{thm}\label{thm4}
Assume $F\in\mathcal{S}$. Then $H\in\mathcal{S}$ if and only if either $D[F]=\emptyset$, or $D[F]\not=\emptyset$ and
(\ref{blxtj}) holds.
\end{thm}

Note that if $G$ has a bounded support then the sufficient and necessary condition above holds trivially. When $G$ has an unbounded support,
there is an analogy between Theorems \ref{1b} and \ref{thm4}. {We feel that under some other premises, the condition that either $D[F]=\emptyset$, or $D[F]\not=\emptyset$ and (\ref{blxtj}) holds, may also be a sufficient and necessary condition for some conclusions, such as $H\in\mathcal{S}$.}
\begin{re}\label{re10}
We make some remarks on Theorem \ref{thm4}:

(i) In general, it is not easy to verify condition (\ref{blxtj}) since $H$ is usually unknown. For this reason, Corollary 1.1 of Tang (2008) gives three useful sufficient conditions (A), (B) and (C) in terms of  the known distributions $F$ and $G$ for (\ref{blxtj}).
We would like to point out that condition (B) that $G\in\mathcal{L}$ may be slightly weakened to
$$G^{*k}\in\mathcal{L} \ \text{for some}\ k\ge2.$$
To see this, recall Theorem 2.1(1a) of Xu et al. (2015), which shows that if $G^{*k}\in\mathcal{L}$ for some $k\ge1$, then, for any constant $d>0$,
$$\overline G(x/d)-\overline G\big((x+1)/d\big)=o\big(\overline {G^{*k}}(x/d)\big).$$
Also notice that
$$\overline {G^{*k}}(x/d)\leq k\overline {G}\big(x/(kd)\big)\le k\overline H(x)/\overline F(kd).$$
Thus, (\ref{blxtj}) holds.
Note that there exist distributions $G$ such that $G\not\in\mathcal{L}$ and  $G^{*k}\in\mathcal{L}$ for some $k\ge2$; see the proof of Theorem 2.2 or Proposition 2.1 in Xu et al. (2015).

(ii) It is known that the subexponentiality of $H$ does not necessarily imply that of $F$ or $G$.
Indeed, Example 2.1 of Tang (2008) shows that the product convolution of two standard exponential distributions is subexponential. This example can be slightly extended to the following:
$$\overline F(x)=\overline G(x)=\textbf{\emph{1}}(x\leq0)+e^{-x^{\alpha}}\textbf{\emph{1}}(x>0) \ \text{for some}\ 1<\alpha<2.$$
Then none of $F$ and $G$ belongs to $\mathcal{L}(\gamma)$ for any $\gamma\ge0$. However, by Lemma 2.1 of Arendarczyk and D\c{e}bicki (2011), we know that
$$\overline H(x)\sim\sqrt{\pi}x^{\alpha/4}e^{-2x^{\alpha/2}}$$
and, therefore, that $H\in\mathcal{S}$. See also Theorem 2.1 of Liu and Tang (2010) for a more general result.

(iii) Similarly, to address the reverse of Theorem \ref{1b}, we will prove two results in Subsection 3.2 below in which
$H\in\mathcal{L}\setminus\mathcal{S}$ but none of $F$ and $G$ needs to belong to the class $\mathcal{L}(\gamma)$ for some $\gamma\ge0$.
\end{re}

Motivated by an open question in Acknowledgment of Tang (2008), now we consider the reverse problem: under what conditions does
the subexponentiality of $H$ imply that of $F$?

\begin{thm}\label{thm5}
Assume $H\in\mathcal{S}$. Then $F\in\mathcal{S}$ if any one of the following
two conditions holds:

(i) $F\in\mathcal{L}$ and
\begin{equation}\label{thm102}
\overline H(x)=O\big(\overline F(x/t)\big)\ \text{for some}\ t\ge1;
\end{equation}

(ii) $G$ is discontinuous and there exists a constant $d\in D[G]$ such that
\begin{equation}\label{thm1002}
\overline H(x)=O\big(\overline F(x/d)\big).
\end{equation}
\end{thm}

\begin{re}\label{re13} (i) %There exist many distributions $F$, $G$ and $H$ satisfying condition (\ref{thm102}) or (\ref{thm1002}). For example, let $F\in\mathcal{R}$ with a regular index $\alpha>0$ and $EY^{\alpha+\varepsilon}<\infty$ for some $\varepsilon>0$, then by the well-known Breiman's (1965) theorem, we have $\overline H(x)\sim EY^{\alpha}\overline F(x)$, which means that both condition (\ref{thm102}) holds for all $t\ge1$ and condition (\ref{thm1002}) holds for all $d\in D[G]$. [[THE DISCUSSION HERE DOESN'T MAKE SENSE TO ME. IF YOU KNOW $F\in\mathcal{R}$ ALREADY, THEN THERE'S NO NEED TO VERIFY $F\in\mathcal{S}$. PLAESE THINK TO SIMPLY DELETE IT.]] More examples can be found
%in Theorem 2.1 of Yang and Wang (2015).
If the distribution $G$ has a finite support $[0,s]$ for some constant $s>0$, then condition (\ref{thm102}) holds for all $t\ge s$.

(ii) One may construct examples in which both  $H$ and $F$ belong to the class $\mathcal{S}$ but none of conditions (\ref{thm102}) and (\ref{thm1002}) holds. Actually, in Example \ref{exam201} in Section 3, for every $t\ge1$ there exists a sequence of positive numbers $\{x_n,n\ge1\}$ such that $x_n\uparrow\infty$ and
\[
\overline{H}(x_{n})/\overline{F}(x_{n}/t)\geq \overline{G}(x_{n})/\overline{F%
}(x_{n}/t)=t^{\alpha }x_{n}\rightarrow \infty ,
\]
as $n\to\infty$, which fails both conditions (\ref{thm102}) and (\ref{thm1002}).

(iii) We may show cases where conditions (\ref{thm102}) and (\ref{thm1002}) are not necessary. For example, let $F$ be a continuous distribution, $F \in \mathcal{L}\setminus\mathcal{S}$, and $G \in \mathcal{S}$. Then by Theorem \ref{thm4}, $H\in\mathcal{S}$, but condition (\ref{thm102}) clearly does not hold. Otherwise, if condition (1.3) holds, then $F\in\mathcal{S}$ by Theorem 1.2, which contradicts to the fact that $F \in \mathcal{L}\setminus\mathcal{S}$. The following result further illustrates this point.
\end{re}

\begin{pron}\label{pron1} Let $F$ be a continuous distribution, either $F\in \mathcal{L}\setminus\mathcal{S}$ or $F \in \mathcal{L}(\gamma)$ for some $\gamma>0$. We can always construct a discontinuous distribution $G$ with unbounded support such that $G\not\in\mathcal{L}$, $H\in\mathcal{S}$, while none of conditions (\ref{thm102}) and (\ref{thm1002}) holds.
\end{pron}

\begin{re}\label{re14}
(i) The proof of Proposition \ref{pron1} to be given in Section 2 below shows a method to construct examples in which $H\in\mathcal{S}$ but none of $F$ and $G$ is subexponential.

(ii) There are actually many examples of distributions belonging to class $\mathcal{L}\setminus\mathcal{S}$; see Pitman(1979), Embrechts and Goldie(1980), Murphree (1989), Leslie(1989), Lin and Wang (2012), Wang et al. (2016), among others.
\end{re}

The rest of this paper consists of three sections. Section 2 gives the proofs of Theorems \ref{thm4},
\ref{thm5} and Proposition \ref{pron1}. Section 3 discusses condition (d) in Theorem 1.A, {and gives two corresponding results for the long-tailed distribution class}. Finally, Section 4 proposes an application of Theorem \ref{thm4} to the asymptotic study of the ruin probability in a discrete-time insurance risk model with stochastic returns.

\section{Proofs}
\setcounter{thm}{0}\setcounter{Corol}{0}\setcounter{lemma}{0}\setcounter{pron}{0}\setcounter{equation}{0}
\setcounter{re}{0}\setcounter{exam}{0}\setcounter{property}{0}\setcounter{defin}{0}
\mbox{}
In this section, we prove Theorems \ref{thm4}, \ref{thm5} and Proposition \ref{pron1}, respectively.
\subsection{Proof of Theorem \ref{thm4}}
The necessity part follows from Theorem \ref{1b}
with $\gamma=0$, since $\mathcal{S}$ is a subset of $\mathcal{L}$. Next, we prove the sufficiency part. Still by Theorem \ref{1b}, we have $H\in\mathcal{L}$.
Then the only statement left to be proved is that: given $F\in\mathcal{S}$ and $H\in\mathcal{L}$, then $H\in\mathcal{S}$.
Following Lemma 2.3 (i) in Cline and Samorodnitsky (1994) and
noting that $\mathcal{L}$ and $\mathcal{S}$ are closed
under scalar multiplication, we only need to prove the result for the case that $Y\ge1$ $a.s.$.

First, note that from Corollary 2.5 of Cline and Samorodnitsky (1994), $F\in\mathcal{S}$ and $H\in\mathcal{L}$, there exist two functions $a_1(\cdot)$ and $a_2(\cdot):[0,\infty)\mapsto(0,\infty)$ satisfying that $a_i(x)\uparrow\infty,\ a_i(x)/x\downarrow 0,\ i=1,2,$
$$\overline F\big(x-a_1(x)\big)\sim \overline F(x)\ \text{and}\ \overline H\big(x-a_2(x)\big)\sim \overline H(x).$$
Let $a(x)=\big(\min\{a_1(x),a_2(x)\}\big)^{1/2}$ for $x\ge0$. From the fact that $a(x)\leq a^2(x)$ for $x$ large enough, we have $a(x)\uparrow\infty,\ a^2(x)/x\downarrow 0,$
$$\overline F(x-a^2(x))\sim \overline F(x)\ \text{and}\ \overline H(x-a(x))\sim \overline H(x).$$

For the given function $a(\cdot)$, {there exist two positive functions $b(\cdot)$ and $c(\cdot)$ supported on $[0,\infty)$ satisfying the following two conditions:
\begin{eqnarray}\label{3002}
c(x)\uparrow\infty,\ b(x)\uparrow\infty,\ c(x)=o\big(b(x)\big)\ \text{and}\ b(x)c(x)=o\big(a(x)\big),
\end{eqnarray}
and
\begin{eqnarray}\label{30002}
\overline G\big(a(x)\big)=o\big(\overline G(2c(x))\big)\ \text{and}\
\overline H\big(a(x)\big)=o\big(\overline F(b(x))\big).
\end{eqnarray}
For brevity, we just construct the function $b(\cdot)$.
%In fact, there exist two positive functions $b_1(\cdot)$ and $c_1(\cdot)$ supported on $[0,\infty)$ such that
%$$\overline G\big(a(x)\big)=o\big(\overline G(2c_1(x))\big)\ \text{and}\ \overline H\big(a(x)\big)=o\big(\overline F(b_1(x))\big).$$
For every $n\ge1$, there exists a positive number $x_n$ such that
$$\overline H\big(a(x)\big)/\overline F(n)<1/n\ \text{for}\ x\ge x_n.$$
Without loss of generality, we set $x_1<x_2<\cdot\cdot\cdot<x_n\uparrow\infty$, and define a function $b_1(\cdot)$ as follows:
$$b_1(x)=\textbf{1}\big(x\in[0,x_{1})\big)+\sum_{n=1}^\infty n\textbf{1}\big(x\in[x_n,x_{n+1})\big)\ \ \text{for}\ x\ge0.$$
Clearly, $b_1(x)\uparrow\infty$ and $\overline H\big(a(x)\big)=o\big(\overline F(b_1(x))\big)$. Further, we define a function: $$b(x)=\min\big\{b_1(x), \big(a(x)\big)^{1/2}\big\}\ \ \text{for}\ x\ge0,$$
which satisfies the above-mentioned requirements that $b(x)\uparrow\infty$ and $\overline H\big(a(x)\big)=o\big(\overline F(b(x))\big)$.}

{Let $(X_i,Y_i)$ for $i=1,2$, be independent copies of $(X,Y)$.} For any $x>0$, we have
\begin{eqnarray}\label{3001}
\overline{H^{*2}}(x)&=&
P\Big(X_1Y_1+X_2Y_2>x,\bigcap_{i=1}^{2}\{1\leq Y_i\leq x/b(x)\}\Big)\nonumber\\
&&\ \ \ \ \ \ \ \  \ \ \ \ \ +P\Big(X_1Y_1+X_2Y_2>x,\bigcup_{i=1}^{2}\{Y_i>x/b(x)\}\Big)\nonumber\\
&= &L_1+L_2.
\end{eqnarray}

We first estimate $L_1$. Clearly, we have
\begin{eqnarray}\label{3003}
L_1&=&P\big(X_1Y_1+X_2Y_2>x,1\leq Y_1\leq Y_2\leq x/b(x)\big)\nonumber\\
&&\ \ \ \ \ \ \ \  \ \ \ \ \ +P\big(X_1Y_1+X_2Y_2>x,1\leq Y_2<Y_1\leq x/b(x)\big)\nonumber\\
&=&L_{11}+L_{12}.
\end{eqnarray}

For $L_{11}$, we have the following decomposition:
\begin{eqnarray}\label{3004}
L_{11}&\leq& P\big(X_1Y_1>x-c(x)a(x),1\leq Y_1\leq Y_2\leq x/b(x)\big)\nonumber\\
& &\ \ \ \ \ \ \ \ +P\big(X_2Y_2>x-c(x)a(x),1\leq Y_1\leq Y_2\leq x/b(x)\big)\nonumber\\
& &\ \ \ \ \ \ \ \ +P\big(X_1Y_1+X_2Y_2>x,c(x)a(x)\leq X_2Y_2\leq x-c(x)a(x),1\leq Y_1\leq Y_2\leq x/b(x)\big)\nonumber\\
&=&L_{111}+L_{112}+L_{113}.
\end{eqnarray}
Note that $a^2(x)/x\downarrow 0$ implies that $a^2(x/z)\geq a^2(x)/z$ for all $z\geq 1$ and $x>0$. Further, by
(\ref{3002}), $\overline F\big((x/z)-a^2(x/z)\big)\sim \overline F\big(x/z\big)$ holds uniformly for $0<z\leq x/b(x)$,
one can hence conclude that
\begin{eqnarray}\label{3005}
L_{111}
&\le&\int_{1}^{x/b(x)}\int_{1}^{y}\overline F\big((x-a^2(x))/z\big)G(dz)G(dy)\nonumber\\
&\leq&\int_{1}^{x/b(x)}\int_{1}^{y}\overline F\big((x/z)-a^2(x/z)\big)G(dz)G(dy)\nonumber\\
&\sim&\int_{1}^{x/b(x)}\int_{1}^{y}\overline F(x/z)G(dz)G(dy)\nonumber\\
&=&P\big(X_1Y_1>x,1\leq Y_1\leq Y_2\leq x/b(x)\big).
\end{eqnarray}
Using a similar approach, we have
\begin{eqnarray}\label{3006}
L_{112}\lesssim P\big(X_2Y_2>x,1\leq Y_1\leq Y_2\leq x/b(x)\big).
\end{eqnarray}
Now, we decompose $L_{113}$ as follows:
\begin{eqnarray}\label{30007}
&&L_{113}\leq P\big(X_1Y_2+X_2Y_2>x,c(x)a(x)<X_2Y_2\leq x-c(x)a(x),c(x)<Y_1\leq Y_2\leq x/b(x)\big)\nonumber\\
&+&P\big(X_1Y_2+X_2Y_2>x,c(x)a(x)<X_2Y_2\leq x-c(x)a(x),1\leq Y_2\leq a(x)\big)\nonumber\\
&+&P\big(X_1Y_1+X_2Y_2>x,c(x)a(x)<X_2Y_2\leq c(x)x/b(x),a(x)< Y_2\leq x/b(x), 1\leq Y_1\leq c(x)\big)\nonumber\\
&+&P\big(X_1Y_2+X_2Y_2>x,\frac{c(x)x}{b(x)}<X_2Y_2\leq x-c(x)a(x),a(x)< Y_2\leq \frac{x}{b(x)}, 1\leq Y_1\leq c(x)\big)\nonumber\\
&&\ \ \ \ \ \ =\sum_{i=1}^{4}L_{113i}.
\end{eqnarray}
Firstly,
according to (\ref{3002}) and $F\in\mathcal{S}$, we have
\begin{eqnarray}\label{30006}
L_{1131}&\leq&P\big((X_1+X_2)Y_2>x, c(x)< Y_2\leq x/b(x), Y_1>c(x)\big)\nonumber\\
&=&\overline{G}\big(c(x)\big)\int_{c(x)}^{x/b(x)}\overline {F^{*2}}(x/y)G(dy)\nonumber\\
&\sim&\overline{G}\big(c(x)\big)\int_{c(x)}^{x/b(x)}2\overline {F}(x/y)G(dy)\nonumber\\
&=&o\big(\overline H(x)\big).
\end{eqnarray}
Secondly, by the facts that $F\in\cal{S}$ and
$$\overline {F^{*2}}(x/y)-\big(1+F(c(x)\big)\overline F(x/y)=o\big(\overline F(x/y)\big)$$
holds uniformly for $1\leq y\leq a(x),$ we have
\begin{eqnarray}\label{3007}
L_{1132}
&=&\int_{1}^{a(x)}\int_{c(x)a(x)/y}^{(x-c(x)a(x))/y}\overline F\big((x/y)-z\big)F(dz)G(dy)\nonumber\\
&\leq&\int_{1}^{a(x)}\Big(\overline {F^{*2}}(x/y)-\overline F(x/y)-\int_{0}^{c(x)a(x)/y}\overline F((x/y)-z)F(dz)\Big)G(dy)\nonumber\\
&\leq&\int_{1}^{a(x)}\Big(\overline {F^{*2}}(x/y)-\big(1+F(c(x)\big)\overline F(x/y)\Big)G(dy)\nonumber\\
&=&o\big(\overline H(x)\big).
\end{eqnarray}
Thirdly, since
$$\overline{H}(x)\geq \overline{F}\big(\frac{x}{c(x)}-\frac{x}{b(x)}\big)\overline{G}\big(\frac{b(x)c(x)}{b(x)-c(x)}\big)$$
and 
$$\frac{b(x)c(x)}{b(x)-c(x)}\leq 2c(x)$$ 
for $x$ large enough, and by (\ref{30002}), we have
\begin{eqnarray}\label{30008}
L_{1133}
&\leq&P\Big(X_1Y_1>x-\frac{c(x)x}{b(x)},a(x)< Y_2\leq \frac{x}{b(x)} ,1\leq Y_1\leq c(x)\Big)\nonumber\\
&\leq&\overline{F}\big(\frac{x}{c(x)}-\frac{x}{b(x)}\big)\overline{G}(a(x))\nonumber\\
&\leq&\overline{H}(x)\overline{F}\big(\frac{x}{c(x)}-\frac{x}{b(x)}\big)\overline{G}\big(a(x)\big)
\Big(\overline{F}\big(\frac{x}{c(x)}-\frac{x}{b(x)}\big)\overline{G}\big(\frac{b(x)c(x)}{b(x)-c(x)}\big)\Big)^{-1}\nonumber\\
&=&o\big(\overline{H}(x)\big).
\end{eqnarray}
Finally, we have
\begin{eqnarray}\label{30009}
L_{1134}
&\leq&\int_{a(x)}^{x/b(x)}\int_{(c(x)x)/(b(x)y)}^{(x-c(x)a(x))/y}\overline F\big((x/y)-z\big)F(dz)G(dy)\nonumber\\
&\leq&\int_{a(x)}^{x/b(x)}\Big(\overline {F^{*2}}(x/y)-\overline F(x/y)-\int_{0}^{(c(x)x)/(b(x)y)}\overline F\big((x/y)-z\big)F(dz)\Big)G(dy)\nonumber\\
&\leq&\int_{a(x)}^{x/b(x)}\Big(\overline {F^{*2}}(x/y)-\big(1+F(c(x))\overline F(x/y)\big)\Big)G(dy)\nonumber\\
&=&o\big(\overline{H}(x)\big).
\end{eqnarray}
Substituting (\ref{30006})-(\ref{30009}) into relation (\ref{30007}) yields that
\begin{eqnarray}\label{30010}
L_{113}&=&o(\overline H(x)).
\end{eqnarray}
Further, substituting (\ref{3005}), (\ref{3006}) and (\ref{30010}) into relation (\ref{3004}) yields that
\begin{eqnarray}\label{3008}
L_{11}&\lesssim&P\big(X_1Y_1>x,1\leq Y_1\leq Y_2\leq x/b(x)\big)+P\big(X_2Y_2>x,1\leq Y_1\leq Y_2\leq x/b(x)\big)\nonumber\\
&&\ \ \ \ \ \ \ \ \ \ +o\big(\overline{H}(x)\big).
\end{eqnarray}

For $L_{12}$, by symmetry, we obtain that
\begin{eqnarray}\label{3009}
L_{12}&\lesssim&P\big(X_1Y_1>x,1\leq Y_2< Y_1\leq x/b(x)\big)+P\big(X_2Y_2>x,1\leq Y_2< Y_1\leq x/b(x)\big)\nonumber\\
&&\ \ \ \ \ \ \ \ \ \ +o\big(\overline{H}(x)\big).
\end{eqnarray}
Clearly, it follows from (\ref{3003}), (\ref{3008}) and (\ref{3009}) that
\begin{eqnarray}\label{30101}
L_1\lesssim 2P\big(X_1Y_1>x,1\leq Y_1\leq x/b(x)\big)+o\big(\overline{H}(x)\big).
\end{eqnarray}

Nexst, we are ready to estimate $L_2$. Clearly, we have
\begin{eqnarray}\label{3010}
L_2&\leq&2P\big(X_1Y_1+X_2Y_2>x,Y_1>x/b(x)\big)\nonumber\\
&\leq&2P\big(X_1Y_1>x-a(x),Y_1>x/b(x)\big)+2P\big(X_2Y_2>a(x),Y_1>x/b(x)\big)\nonumber\\
&=&2L_{21}+2L_{22}.
\end{eqnarray}

For $L_{21}$, it follows that
\begin{eqnarray}\label{3011}
L_{21}&=&P\big(X_1Y_1>x-a(x))-P(X_1Y_1>x-a(x),Y_1\leq x/b(x)\big)\nonumber\\
&\leq&P\big(X_1Y_1>x-a(x)\big)-P\big(X_1Y_1>x,Y_1\leq x/b(x)\big)\nonumber\\
&=&\overline H(x)+o\big(\overline H(x)\big)-P\big(X_1Y_1>x,Y_1\leq x/b(x)\big)\nonumber\\
&=&P\big(X_1Y_1>x,Y_1>x/b(x)\big)+o\big(\overline H(x)\big).
\end{eqnarray}

And for $L_{22}$, from (\ref{30002}), we have
\begin{eqnarray}\label{3012}
L_{22}&=&\overline H(x)\overline H(a(x))\overline G(x/b(x))/\overline H(x)\nonumber\\
&\leq&\overline H(x)\overline H(a(x))/\overline F(b(x))\nonumber\\
&=&o\big(\overline{H}(x)\big).
\end{eqnarray}

Hence by (\ref{3001}) and (\ref{30101})-(\ref{3012}), it holds that
\begin{eqnarray*}
\overline{H^{*2}}(x)\lesssim 2\overline{H}(x),
\end{eqnarray*}
which, combining with Theorem 2.11 in Foss et al. (2013), implies that $H\in\cal{S}$.

\subsection{Proof of Theorem \ref{thm5}}
(i) For some $t\ge1$ in (1.3), let $F_1$ be a distribution such that $\overline{F_1}(x)=\overline{F}(x/t)$ for all $x$. By $H\in\mathcal{S}$, $F\in\mathcal{L}$ and (\ref{thm102}), we have $F_1\in\mathcal{L}$ and $\overline{F_1}(x)\asymp\overline{H}(x)$, thus by Theorem 3.11 in Foss et al. (2013), $F_1\in\mathcal{S}$. Therefore,
$$\overline{F^{*2}}(x)=\overline{F_1^{*2}}(tx)\sim2\overline{F_1}(tx)=2\overline{F}(x),$$
namely $F\in\mathcal{S}$.

(ii) From the following fact:
\begin{eqnarray*}
\overline{H}(x-2)-\overline{H}(x)&\ge&\int_{(d-(d/x),d]}\big(\overline F((x-2)/y)-\overline F(x/y)\big)G(dy)\nonumber\\
&\ge&\big(\overline F((x-2)/(d-(d/x)))-\overline F(x/d)\big)G\{d\}\nonumber\\
&\ge&\big(\overline F((x-1)/d)-\overline F(x/d)\big)G\{d\},
\end{eqnarray*}
$H\in\mathcal{S}$ and (\ref{thm1002}), we have
\begin{eqnarray*}
\overline F\big((x-1)/d\big)-\overline F(x/d)=o\big(\overline H(x)\big)=o\big(\overline F(x/d)\big),
\end{eqnarray*}
which implies $F\in\mathcal{L}$.
Then by (i),
we can get $F\in\mathcal{S}$.

\subsection{Proof of Proposition \ref{pron1}} Let $X$ and $Y_0$
be two independent nonnegative random variables with continuous distributions $F$ and $G_1$, such that $G_1\in\mathcal{S}$ and either $F\in\mathcal{L}(\gamma)$ for some $\gamma>0$ or $F\in \mathcal{L}\setminus\mathcal{S}$.  Denote the distribution of $XY_0$ by $H_1$. Since $G_1$ is continuous, by Theorem \ref{thm4}, we have $H_1\in\mathcal{S}$.

Using the method in Example 3.1 of Xu et al. (2016), we can construct a new random variable $Y$ with discontinuous distribution $G$, satisfying $G\not\in\mathcal{L}$ and $\overline {G}(x)\asymp\overline {G_1}(x)$, which implies $\overline {H}(x)\asymp\overline {H_1}(x)$. Further, we can get $H\in\mathcal{L}$ directly from Theorem \ref{1b}. Thus we know that $H\in\mathcal{S}$.

Clearly, conditions (\ref{thm102}) and (\ref{thm1002}) are not satisfied in this case.

\section{Discussion}
\setcounter{thm}{0}\setcounter{Corol}{0}\setcounter{lemma}{0}\setcounter{pron}{0}\setcounter{equation}{0}
\setcounter{re}{0}\setcounter{exam}{0}\setcounter{property}{0}\setcounter{defin}{0}

Firstly, we further discuss condition (d) in Theorem 1.A by providing (1) an example and (2) two sufficient conditions for $H\in \mathcal{L}$, respectively.

\subsection{On condition (d)}

The following example shows that, when $F\in\mathcal{S}$, condition (d) is unnecessary for $H\in\mathcal{S}$.
\begin{exam}\label{exam201}
Choose any constants $\alpha>0$ and $x_1>4^{\alpha}$. For all integers $n\ge1$, let $x_{n+1}=x_n^{1+1/\alpha}$. Clearly, $x_{n+1}>4x_n$ and $x_n\to\infty$, as $n\to\infty$. We define two distributions $F$ and $G$ such that
\begin{eqnarray*}
&&\overline G(x)=\textbf{\emph{1}}(x<0)+
\big(1+(x_1^{-\alpha-1}-x_1^{-1})x\big)\textbf{\emph{1}}(0\leq x<x_1)\nonumber\\
&&+\sum_{n=1}^{\infty}\Big(\big(x_n^{-\alpha}+
(x_n^{-\alpha-2}-x_n^{-\alpha-1})(x-x_n)\big)\textbf{\emph{1}}(x_n\leq x<2x_n)+x_n^{-\alpha-1}\textbf{\emph{1}}(2x_n\leq x<x_{n+1})\Big)
\end{eqnarray*}
and
$$\overline F(x)=\textbf{\emph{1}}(x<1)+\textbf{\emph{1}}(x\ge1)x^{-\alpha-1}.$$
\end{exam}
Because the distribution $F\in\mathcal{R}\subset\mathcal{S}$ and is continuous, by Theorem \ref{thm4}, $H\in\mathcal{S}$.
From
\begin{eqnarray*}
\overline G(2x_n-1)/ \overline G(2x_n)=2-x_n^{-1}\to2,
\end{eqnarray*}
as $n\to\infty,$ we know $G\notin\mathcal{L}$. Direct calculation leads to
\begin{eqnarray*}
\overline H(x)=x^{-\alpha-1}\int_{0}^{x}y^{\alpha+1}G(dy)+\overline G(x).
\end{eqnarray*}
Thus, when $x\in[x_n,2x_n)$,
\begin{eqnarray}\label{26}
\overline H(x)&=&x^{-\alpha-1}\Big(\int_{0}^{x_1}(x_1^{-1}-x_1^{-\alpha-1})
y^{\alpha+1}dy+\sum_{i=1}^{n-1}\int_{x_i}^{2x_i}(x_i^{-\alpha-1}-x_i^{-\alpha-2})
y^{\alpha+1}dy\nonumber\\
& &\ \ \ \ \ \ \ +\int_{x_n}^{x}(x_n^{-\alpha-1}-x_n^{-\alpha-2})
y^{\alpha+1}dy\Big)+\overline G(x);
\end{eqnarray}
when $x\in[2x_n,x_{n+1})$,
\begin{eqnarray}\label{27}
\overline H(x)&=&x^{-\alpha-1}\Big(\int_{0}^{x_1}(x_1^{-1}-x_1^{-\alpha-1})
y^{\alpha+1}dy+\sum_{i=1}^{n}\int_{x_i}^{2x_i}(x_i^{-\alpha-1}-x_i^{-\alpha-2})
y^{\alpha+1}dy\Big)+\overline G(x).
\end{eqnarray}

Now, we prove $H\in\mathcal{D}$. In fact,
when $x\in[x_n,4x_n)$, by (\ref{27}), we have
\begin{eqnarray*}
\overline H(x/2)/\overline H(x)&\leq&\overline H(x_n/2)/\overline H(4x_n)\nonumber\\
&\leq&2^{3+3\alpha}+2^{2+2\alpha}\nonumber\\
&<&\infty;
\end{eqnarray*}
when $x\in[4x_n,x_{n+1})$, by the inequality
$$(a+b)/(c+d)\leq (a/b)+(c/d)$$
for all $a,b,c,d>0$ and (\ref{27}), we know that
\begin{eqnarray*}
\overline H(x/2)/\overline H(x)&\leq&2^{\alpha+1}+1<\infty.
\end{eqnarray*}

From (\ref{26}), we know that $\overline G(x_n)/\overline H(x_n)\to1,$ as $n\to\infty$. Thus, for any $b>0$, $\overline G(x)=o(\overline H(bx))$ doesn't hold due to $H\in\mathcal{D}$, which means that for any positive function $a(x)$ satisfying $a(x)\uparrow\infty$ and $a(x)/x\downarrow 0$, condition (d) does not hold.

\subsection{On the long-tailed distributions}

Next, we gives the following two results, which show that even if two distributions $F$ and $G$ may not belong to the class $\mathcal{L}(\gamma)$ for any $\gamma\ge0$, the product convolution $H$ is still long-tailed under certain conditions.
\begin{thm}\label{thm7}
Let $F$ and $G$ be two continuous distributions. If there exists a function $a(\cdot):[0,\infty)\mapsto(0,\infty)$ satisfying
$a(x)\uparrow\infty$, $x/a(x)\uparrow\infty$,
\begin{equation}\label{cw01}\overline G\big(a(x)\big)=O\big(\overline H(x)\big)
\end{equation}
and
\begin{equation}\label{cw011}
\overline F\big(x/a(x)\big)=O\big(\overline H(x)\big),
\end{equation}
then $H\in\mathcal{L}$.
\end{thm}
\proof To prove the result, we estimates the following tail distribution function. For $x>1$,
\begin{eqnarray}
&&\overline H(x-1)=P\big(XY>x-1,Y\leq a(x)\big)+P\big(XY>x-1,Y> a(x)\big)\nonumber\\
&\leq&P\big(XY>x-1,X\geq (x-1)/a(x)\big)+P\big(XY>x-1,Y> a(x)\big)\nonumber\\
&=&\Big(\int_{(x-1)/a(x)}^{(x-1)/a(x-1)}+\int_{(x-1)/a(x-1)}^{x/a(x)}+\int_{x/a(x)}^{\infty}\Big)\overline G\big((x-1)/y\big)F(dy)+\int_{a(x)}^{\infty}\overline F\big((x-1)/y\big)G(dy)\nonumber\\
&\leq&\int_{x/a(x)}^{\infty}\overline G\big((x-1)/y\big)F(dy)+\int_{a(x)}^{\infty}\overline F\big((x-1)/y\big)G(dy)\nonumber\\
& &\ \ \ \ \ \ \ \ +\overline F\big((x-1)/a(x)\big)\overline G\big(a(x-1)\big)+\overline F\big((x-1)/a(x-1)\big)\overline G\big((x-1)a(x)/x\big).\label{cw02}
\end{eqnarray}
Clearly, it follows from (\ref{cw01}) and (\ref{cw011}) that
\begin{equation}\label{cw03}
\overline F\big((x-1)/a(x)\big)\overline G\big(a(x-1)\big)+\overline F\big((x-1)/a(x-1)\big)\overline G\big((x-1)a(x)/x\big)=o\big(\overline{H}(x-1)\big).
\end{equation}
Moreover, since $F$ and $G$ are continuous on $[0,\infty)$, they are uniformly continuous on $[0,\infty)$. Hence, for any fixed $\epsilon>0$ and for $x$ large enough, we have that
$$\overline G\big((x-1)/y\big)<\overline G(x/y)+\epsilon$$
and
$$\overline F\big((x-1)/y\big)<\overline F(x/y)+\epsilon$$
uniformly holds for $y>x/a(x)$. Combining with (\ref{cw01}), it follows that
\begin{eqnarray}
 &&\int_{x/a(x)}^{\infty}\overline G\big((x-1)/y\big)F(dy)+\int_{a(x)}^{\infty}\overline F\big((x-1)/y\big)G(dy)\nonumber\\
 &=& \int_{x/a(x)}^{\infty}\overline G(x/y)F(dy)+\int_{a(x)}^{\infty}\overline F(x/y)G(dy)+o(\overline{H}(x))\nonumber\\
 &=&\overline H(x)+o(\overline H(x)).\label{cw04}
\end{eqnarray}
Thus $H\in\mathcal{L}$ follows from (\ref{cw02})-(\ref{cw04}).
\hfill$\Box$

\begin{thm}\label{thm8}
Let $F$ and $G$ be two continuous distributions. If
\begin{eqnarray}\label{801}
\overline F(x-1/x)\sim\overline F(x) \ \text{and} \ \overline G(x-1/x)\sim \overline G(x)
\end{eqnarray}
and there exists a function $a(\cdot):[0,\infty)\mapsto(0,\infty)$ satisfying
$a(x)\uparrow\infty$, $x/a(x)\uparrow\infty$,
\begin{equation}\label{cw11}
\overline G\big(a(x)\big)=O\big(\overline H(x)\big),
\end{equation}
and
\begin{equation}\label{cw111}
\overline F\big(a(x)\big)=O\big(\overline H(x)\big),
\end{equation}
then $H\in\mathcal{L}$.
\end{thm}
\proof~~Without loss of generality, we assume $a(x)>\sqrt{x}$. Consider the following expression %equality
\begin{eqnarray}\label{cw13}
&&\overline H(x-1)=P(XY>x-1,Y\leq \sqrt{x})+P(XY>x-1,Y> \sqrt{x})\nonumber\\
&=& P(XY>x-1,X>\sqrt{x})+P(XY>x-1,Y>\sqrt{x})\nonumber\\
& &\ \ \ \ +P\big(XY>x-1,(x-1)/\sqrt{x}<X\leq \sqrt{x}\big)-P\big(X>(x-1)/\sqrt{x},Y>\sqrt{x}\big)\nonumber\\
&\leq&\Big(\int_{\sqrt{x}}^{a(x)}+\int_{a(x)}^{\infty}\Big)\overline G((x-1)/y)F(dy)+\Big(\int_{\sqrt{x}}^{a(x)}+\int_{a(x)}^{\infty}\Big)\overline F((x-1)/y)G(dy)\nonumber\\
& &\ \ \ \ +\big(\overline F\big((x-1)/\sqrt{x}\big)-\overline F(\sqrt{x})\big)\big(\overline G\big((x-1)/\sqrt{x}\big)-\overline G(\sqrt{x})\big)-\overline F(\sqrt{x})\overline G(\sqrt{x}).
\end{eqnarray}
From (\ref{801})  and $\overline H(x)\geq \overline F(\sqrt{x})\overline G(\sqrt{x})$, it is easy to know that
\begin{eqnarray}\label{cw15}
\big(\overline F((x-1)/\sqrt{x})-\overline F(\sqrt{x})\big)\big(\overline G((x-1)/\sqrt{x})-\overline G(\sqrt{x})\big)=o\big(\overline H(x)\big).
\end{eqnarray}
Then by
(\ref{801}), for sufficiently large $x$, we have
\begin{eqnarray}\label{cw16}
&&\int_{\sqrt{x}}^{a(x)}\overline G((x-1)/y)F(dy)+\int_{\sqrt{x}}^{a(x)}\overline F((x-1)/y)G(dy)\nonumber\\
&\leq&\int_{\sqrt{x}}^{a(x)}\overline G\Big(\frac{x}{y}-\frac{y}{x}\Big)F(dy)+\int_{\sqrt{x}}^{a(x)}\overline F \Big(\frac{x}{y}-\frac{y}{x}\Big)G(dy)\nonumber\\
&=&\int_{\sqrt{x}}^{a(x)}\overline G(x/y)F(dy)+\int_{\sqrt{x}}^{a(x)}\overline F(x/y)G(dy)+o(\overline{H}(x)).
\end{eqnarray}
Since $F$ and $G$ are continuous on $[0,\infty)$, then they are uniformly continuous on $[0,\infty)$. Using a similar approach as
in Theorem \ref{thm7} and combining with (\ref{cw11}) and (\ref{cw111}), it follows that
\begin{eqnarray}\label{cw17}
&&\int_{a(x)}^{\infty}\overline G\big((x-1)/y\big)F(dy)+\int_{a(x)}^{\infty}\overline F\big((x-1)/y\big)G(dy)\nonumber\\
&=&\int_{a(x)}^{\infty}\overline G(x/y)F(dy)+\int_{a(x)}^{\infty}\overline F(x/y)G(dy)+o\big(\overline{H}(x)\big).
\end{eqnarray}
Plugging (\ref{cw15})-(\ref{cw17}) into relation (\ref{cw13}) yields that
\begin{eqnarray*}
\overline H(x-1)&=&\int_{\sqrt{x}}^{\infty}\overline G(x/y)F(dy)+\int_{\sqrt{x}}^{\infty}\overline F(x/y)G(dy)-\overline F(\sqrt{x})\overline G(\sqrt{x})+o(\overline H(x))\nonumber\\
&=&\overline H(x)+o\big(\overline H(x)\big),
\end{eqnarray*}
which implies that $H\in\mathcal{L}$.\hfill$\Box$

\begin{re}\label{re31}
We takes two distributions $F$ and $G$ as follows:
$$\overline F(x)=\overline G(x)=\textbf{\emph{1}}(x\leq0)+e^{-x^2}\textbf{\emph{1}}(x>0),$$
for all $x$. Then bu Lemma 2.1 of Arendarczyk and D\c{e}bicki (2011), we have
$$\overline H(x)\sim\sqrt{\pi x}e^{-2x},$$
which yields that $H\in\mathcal{L}(2)$. Notice that (\ref{cw11}) holds for all $a(x)=x^\beta$, where $\beta\in(1/2, 1)$, but (\ref{801}) does not hold since 
$$\overline F(x-1/x)\sim e^{2}\overline F(x).$$ %Then $H\in\mathcal{L}(2)$.
This shows that condition (\ref{801}) in Theorem \ref{thm8} may be necessary for $H\in\mathcal{L}$ in a certain sense.
\end{re}

\section{Application
}
\setcounter{thm}{0}\setcounter{Corol}{0}\setcounter{lemma}{0}\setcounter{pron}{0}\setcounter{equation}{0}
\setcounter{re}{0}\setcounter{exam}{0}\setcounter{property}{0}\setcounter{defin}{0}

In this section, we consider a discrete-time insurance risk model with both insurance and financial
risks, which was proposed by Nyrhinen (1999, 2001).

For each $i\ge1$, within period $i$ we denote the net insurance loss (i.e. the total claim amount minus the total premium income) by a real-valued r.v. $Z_i$ with common distribution $F_0$ supported on $(-\infty,\infty)$. Further, we denote the distribution of $X=Z_1^+=Z_1\textbf{1}(Z_1\ge0)$ by $F$. Suppose that the insurer makes both risk-free and risky investments, which lead to an overall stochastic discount factor
$Y_i$ with common distribution $G$ supported on $(0,\infty)$ or $(0,s]$ for some $0<s<\infty$, from time $i$ to time $i-1,\ i\ge1$. We assume the $\{Z_i,i\ge1\}$ and $\{Y_i,i\ge1\}$ are two sequences of independent random variables, and $\{Z_i,i\ge1\}$ is independent of $\{Y_i,i\ge1\}$. They are called insurance risks and financial risks respectively in Tang and Tsitsiashvili (2003, 2004). For each positive integer $n$, the sum
\begin{eqnarray*}
S_n=\sum_{i=1}^{n}Z_i\prod_{j=1}^{i}Y_j
\end{eqnarray*}
represents the stochastic discount value of aggregate net losses up to time $n$.
Then the finite-time ruin probability by time $n$ and the infinite-time ruin probability respectively can be defined as
\begin{eqnarray*}
\psi(x;n)=P\big(\max_{1\le k \le n}S_k>x\big)
\end{eqnarray*}
and
\begin{eqnarray*}
\psi(x)=P\big(\sup_{k\ge1}S_k>x\big),
\end{eqnarray*}
where $x\ge0$ is interpreted as the initial wealth of the insurer.

There are a lot of papers studying the asymptotic behavior of the ruin probability in this model. The reader can refer to Goovaerts et al. (2005), Zhang et al. (2009), Chen (2011), Yang and Wang (2013), Huang et al. (2014), Li and Tang (2015), etc. By using Theorem \ref{thm4} of this paper, we can get a different version of the finite-time ruin probability which is stated later. To this end, we denote the distribution of $Z_{i}^+\prod_{j=1}^{i}Y_j$ by $H_i,\ i\ge1$. Clearly, $H_1=H$.

In the risk model, assume that $F\in\mathcal{S}$ and $G$ are supported on $(0,\infty)$ or $[0,s]$ for some $0<s<\infty$, then $S_n$ follows a subexponential distribution for each $n\in\mathbb{N}$ if and only if $D[F]=\emptyset$, or $D[F]\not=\emptyset$ and %it
(\ref{blxtj}) holds.
Further,
\begin{eqnarray*}\label{5001}
\psi(x;n)\sim P(S_n>x)\sim \sum_{i=1}^{n}\overline H_i(x)
\end{eqnarray*}
if condition (\ref{blxtj2}) is satisfied.

Recall that condition (\ref{blxtj2}) holds if and only if there exists a function $a(\cdot):[0,\infty)\mapsto(0,\infty)$ such that
\begin{eqnarray}\label{5000}
a(x)\uparrow\infty,\ a(x)/x\downarrow 0\ \mbox{and} \ \overline G\big(a(x)\big)=o\big(\overline H(x)\big).
\end{eqnarray}
Tang (2006) shows that there exist two distributions $F$ and $G$, as well as a positive function $a(\cdot)$, such that condition (\ref{5000}) holds, while the condition fails to hold for any function satisfying condition (c) in Theorem \ref{1a}. Therefore, the conditions of the above-mentioned result are more relaxed than the existing ones. However, we omit the proof of the result, because the method of the proof is standard.

As for the infinite time ruin probability in the risk model, there are only a few accurate results, see, for example, Theorem 4.2 of Yang and Wang (2013) for $F\in\mathcal{R}$ with some indicator $\alpha>0$. For the case that $F\in\mathcal{S}\setminus\mathcal{R}$, we do not find any relevant results. Here we can provide a result for the asymptotic lower bound of the infinite time ruin probability.

\begin{pron}\label{pron001}
If $F\in\mathcal{A}$ and $G$ is supported on $(0,1]$, then
\begin{eqnarray}\label{pron01}
\liminf\psi(x)\Big/\sum_{i=1}^{\infty}\overline H_i(x)\ge1.
\end{eqnarray}
\end{pron}
\proof Since $F\in\mathcal{A}$, for any fixed $\lambda>1$,
$$a=a(F,\lambda)=\limsup \overline F(\lambda x)/\overline F(x)<1.$$
So for any $0<\varepsilon<1-a$, there is a constant $x_0=x_0(F,\lambda)$ such that for all $x\ge x_0$,
$$\overline F(\lambda x)\leq (a+\varepsilon)\overline F(x).$$
Then it holds uniformly for all $i\ge1$ and $x\ge x_0$ that
\begin{eqnarray}\label{pron00004}
\overline H_i(\lambda x)\big/\overline H_i(x)&=&\int_{0}^{1}\overline F(\lambda x/y)P\Big(\prod_{j=1}^iY_j\in dy\Big)\Big/\int_{0}^{1}\overline F(x/y)P\Big(\prod_{j=1}^iY_j\in dy\Big)\nonumber\\
&\leq&\sup_{0<y\leq 1}\overline F(\lambda x/y)/\overline F(x/y)\nonumber\\
&\leq& a+\varepsilon.
\end{eqnarray}
Thus
\begin{eqnarray}\label{pron00005}
a_i=a_i(H_i,\lambda)=\limsup\overline H_i(\lambda x)/\overline H_i(x)=a\ \ \text{for all}\ \ i\ge1.
\end{eqnarray}
Further, we denote $p=p(G,\lambda)=G((0,1/\lambda])$ and $q=1-p$, then by (\ref{pron00004}) and (\ref{pron00005}), for all $i\ge1$ and $x\ge x_0$, we have
\begin{eqnarray}\label{pron006}
\overline H_{i+1}(x)&=&\Big(\int_{0}^{1/\lambda}+\int_{1/\lambda}^{1}\Big)\overline H_{i}(x/y)G(dy)\nonumber\\
&\leq& p\overline H_{i}(\lambda x)+q\overline H_{i}(x)\nonumber\\
&\leq& (pa+p\varepsilon+q)\overline H_{i}(x)\nonumber\\
&\leq& (pa+p\varepsilon+q)^{i}\overline H(x).
\end{eqnarray}
By (\ref{pron006}) and the fact that $0<pa+p\varepsilon+q<1$, we have
\begin{eqnarray*}
\liminf\psi(x)\Big/\sum_{i=1}^{\infty}\overline H_i(x)&\geq&\liminf\psi(x,n)\Big(1-\limsup\sum_{i=n+1}^{\infty}\overline H_i(x)/\overline H(x)\Big)\Big/\sum_{i=1}^{n}\overline H_i(x)\nonumber\\
&\geq&1-\sum_{i=n+1}^{\infty}(pa+p\varepsilon+q)^{i-1}\nonumber\\
&\to&1,
\end{eqnarray*}
as $n\to\infty.$ Then (\ref{pron01}) is proved.\hfill$\Box$

In general, when $F\in\mathcal{S}\setminus\mathcal{R}$, it is a very interesting and difficult task to obtain asymptotic upper bounds of the infinite time ruin probabilities. Even in a very simple case, the series $\sum_{i=1}^{\infty}\overline H_i(x)$ is divergent. For example,
if $G$ is supported on $[s_1,s_2]$ or $[s_1,\infty)$ for some $1\leq s_1<s_2<\infty$, then for any distribution $F$ supported on $[0,\infty)$ and any integer $n\ge1$,
\begin{eqnarray*}
\sum_{i=1}^{\infty}\overline H_i(x)\geq\sum_{i=1}^{n}\overline F(x/s_1^{i})
\ge n\overline F(x)\to\infty,
\end{eqnarray*}
as $n\to\infty$.

\vspace{1cm}
\noindent{\bf Acknowledgements.} The authors are grateful to Professor Gennady Samorodnitsky for his careful reading of the proof of Theorem \ref{thm4} in the original version of this paper with valuable comments and suggestions. The authors are also grateful to Professor Sergey Foss and Professor Qihe Tang for their detailed and valuable comments and suggestions for the whole content of this paper. Finally, the authors are grateful to two referees for their valuable comments and suggestions.

\end{document}